\documentclass{llncs}
\usepackage{latexsym}
\def\Hide#1{\relax}
\def\Func#1{{\mathsf{#1}}}
\def\SC#1{\Func{V}(#1)}
\def\PZ#1{\Func{Z_{+}}(#1)}

\def\X{\mathsf{X}}
\def\Arr{\rightarrow}

\newtheorem{Theorem}{Theorem}
\newtheorem{Lemma}[Theorem]{Lemma}

\def\Proof{\par \noindent{\bf Proof: }}

\DeclareSymbolFont{AMSa}{U}{msa}{m}{n}
\DeclareMathSymbol\therefore{\mathrel}{AMSa}{"29}

\title{Descartes' Rule of Signs by an Easy Induction}
\author{R.D. Arthan}
\institute{Lemma 1 Ltd. 2nd Floor, 31A Chain Street, Reading UK.  RG1 2HX\\
\email{rda@lemma-one.com}}
\makeindex
\begin{document}
\maketitle
\begin{abstract}
If $c$ is a positive number, Descartes' rule of signs implies that multiplying a polynomial $f(x)$ by $c - x$ introduces an odd number of changes of sign in the coefficients.
We turn this around, proving this fact about sign changes inductively and deriving Descartes' rule from it.
\end{abstract}

If $a_0, \ldots, a_n$ is a list of real numbers, the number of {\em sign changes} in $a_0, \ldots, a_n$, written $\SC{a_0, \ldots, a_n}$, is the number of pairs of indices $(p, q$) such that $p < q$, $a_pa_q < 0$ and $a_i = 0$ for $p < i < q$.
If $f(x) = a_0 + a_1x + a_2 x^2 + \ldots + a_nx^n$ is a polynomial with real coefficients write $\SC{f}$ for $\SC{a_0, \ldots, a_n}$. If $f(x)$ is a non-zero polynomial, write $\PZ{f}$ for the number of positive roots of $f(x)$ counted according to their multiplicity, i.e., $\PZ{f}$ is the sum over all positive real numbers, $c$, such that $f(c) = 0$, of the largest integer, $m(c)$, such that $(c - x)^{m(c)}$ divides $f(x)$.
Descartes' rule of signs is the following theorem:
\begin{Theorem}\label{TheTheorem}
If $f$ is a non-zero polynomial, $\SC{f} - \PZ{f}$ is even and non-negative.
\end{Theorem}
If $\SC{f}$ is odd, one can write $f(x) = x^mg(x)$, where $g(x)=a_0 + a_1 x + \ldots + a_nx^n$ and $a_0a_n < 0$, but then $g(0)g(x)/x^n$ tends to $a_0a_n$ as $x$ tends to infinity, so that $g(0)$ and $g(x)$ have opposite signs for all sufficiently large $x$, and so, by the intermediate value theorem, $g$ and hence $f$ have at least one positive root.
Thus if $\PZ{f} = 0$, $\SC{f}$ must be even.
By induction on the number of distinct roots, the theorem therefore reduces to the following lemma:
\begin{Lemma}\label{TheLemma}
Let $f$ and $g$ be non-zero polynomials such that $f(x) =(c-x)^mg(x)$ for some positive real number $c$ and positive integer $m$, then $\SC{f} - \SC{g} - m$  is even and non-negative.
\end{Lemma}
\Proof
I claim that if $f(x) = (c - x)g(x)$, then $\SC{f} - \SC{g}$ is positive and odd, from which the lemma follows by induction on $m$.
Let $f(x) = a_0 + a_1x + \ldots + a_nx^n$ and $g(x) = b_0 + b_1x + \ldots + b_{n-1}x^{n-1}$.
Since $\SC{f}$ and $\SC{g}$ are unchanged if we multiply the coefficients $a_i$ and $b_i$ by $c^i$, we may assume that $c = 1$.

Now $a_0 + a_1x + \ldots + a_nx^n = (1-x)(b_0 + b_1x + \ldots + b_{n-1}x^{n-1})$ iff the following equations hold:
\begin{eqnarray*}
b_0 &=& a_0 \\
b_1 &=& a_0 + a_1 \\
b_2 &=& a_0 + a_1 + a_2 \\
    &\vdots& \\
b_{n-1} &=& a_0 + a_1 + a_2 + \ldots + a_{n-1}\\
0 &=& a_0 + a_1 + a_2 + \ldots + a_{n-1} + a_n
\end{eqnarray*}
To complete the proof of the lemma, let us show, by induction on $n \ge 1$, that, if $a_0, \ldots, a_n$ and $b_0, \ldots, b_{n-1}$  are any lists of real numbers with $a_n \not= 0$ satisfying the above equations, then $\SC{a_0, \ldots, a_n} - \SC{b_0, \ldots, b_{n-1}}$ is positive and odd.

If $n = 1$, then $b_0 = a_0 = -a_1 \not= 0$, whence $\SC{a_0, a_1} = 1$ and $\SC{b_0} = 0$, so that $\SC{a_0, a_1} - \SC{b_0} = 1$ which is certainly positive and odd.

Now let $a_0, \ldots, a_n$ and $b_0, \ldots, b_{n-1}$ satisfy the above equations for some $n \ge 2$.
Ignoring the first equation, we have:
\begin{eqnarray*}
b_1 &=& (a_0 + a_1) \\
b_2 &=& (a_0 + a_1) + a_2 \\
    &\vdots& \\
b_{n-1} &=& (a_0 +a_1) + a_2 + \ldots + a_{n-1}\\
0 &=& (a_0 + a_1) + a_2 + \ldots + a_{n-1} + a_n,
\end{eqnarray*}
Since $n \ge 2$ and $a_n \not= 0$, the inductive hypothesis applies to the above equations and we have that
$\SC{a_0 + a_1, a_2, \ldots, a_n} - \SC{b_1, \ldots, b_{n-1}}$ is positive and odd.
Let $\alpha = \SC{a_0, a_1, \ldots, a_n} - \SC{a_0 + a_1, a_2, \ldots, a_n}$ and $\beta = \SC{b_0, b_1, \ldots, b_{n-1}} - \SC{b_1, \ldots, b_{n-1}}$.
We need to show that  $\SC{a_0, a_1, \ldots, a_n} - \SC{b_0, b_1, \ldots, b_{n-1}}$ is also positive and odd, or equivalently, that $\alpha - \beta$ is even and non-negative.

As $a_n \not=0$, there is a least $p$ such that $p > 1$ and $a_p \not = 0$.
If, $b_1 = a_0 + a_1 = 0$, then also $p$ is the least $p \ge 1$ such that $b_p \not= 0$.
So we may now calculate $\alpha$ and $\beta$ for each combination of the signs ($+$, $0$, or $-$) of the numbers $a_0 = b_0$, $a_1$, $b_0 = a_0 + a_1$ and $a_p$.
Now, by construction, $a_p \not=0$, and, multiplying the $a_i$ and $b_i$ by $-1$ if necessary, we may assume $a_0 \ge 0$.
The remaining cases are shown in the following table.
\[
\begin{array}{cc|c|c|c||ccc|ccc||c|c}
                    & a_0 & a_1 & b_1 & a_p &  \multicolumn{3}{c|}{\mathbf{A}}  &    \multicolumn{3}{c||}{\mathbf{B}}   & \alpha & \beta \\\hline
\mbox{{\bf (i)}}   &  0  & \X  & \X  & \X  &  PQ &\Arr& 0PQ   &  RS &\Arr& 0RS   &    0   &   0   \\\hline
\mbox{{\bf (ii)}}  & \X  &  0  & \X  & \X  &  PQ &\Arr& P0Q   &  RS &\Arr& RRS   &    0   &   0   \\\hline
\mbox{{\bf (iii)}} &  +  &  +  & (+) & \X  &  +Q &\Arr& +{+}Q &  +S &\Arr& +{+}S &    0   &   0   \\\hline
\mbox{{\bf (iv)}}  &  +  &  -  &  0  &  +  &  0+ &\Arr& +{-}+ &  0+ &\Arr& +0+   &    2   &   0   \\\hline
\mbox{{\bf (v)}}   &  +  &  -  &  0  &  -  &  0- &\Arr& +{-}- &  0- &\Arr& +0-   &    1   &   1   \\\hline
\mbox{{\bf (vi)}}  &  +  &  -  &  +  &  +  &  ++ &\Arr& +{-}+ &  +S &\Arr& +{+}S &    2   &   0   \\\hline
\mbox{{\bf (vii)}} &  +  &  -  &  +  &  -  &  +- &\Arr& +{-}- &  +S &\Arr& +{+}S &    0   &   0   \\\hline
\mbox{{\bf (viii)}}&  +  &  -  &  -  & \X  &  -Q &\Arr& +{-}Q &  -S &\Arr& +{-}S &    1   &   1
\end{array}
\]

In this table the first four columns give a condition on the signs of the quantities $a_0$, $a_1$, $b_0$ and $a_p$
In these columns, $\X$ means that any sign will do and a sign in brackets is one that is forced by the conditions on the other quantities.

The columns labelled $\mathbf{A}$ and $\mathbf{B}$ show how the sign changes in the sequences of $a_i$ and $b_i$ are transformed in the passage from $a_0 + a_1, \ldots, a_p$ to $a_0, a_1, \ldots, a_p$ or from $b_1, \ldots, b_q$ to $b_0, b_1, \ldots, b_q$, where $q$ is the least $q > 1$ such that $b_q \not= 0$ (so that, as observed above, $q = p$ and $b_q = a_p$ when $b_1 = 0$).
In these two columns $P$, $Q$, $R$ and $S$ are used as variables ranging over the set $\{+, 0, -\}$ in cases where the sign in question is not determined by the sign condition determined by the first four columns.

So, for example, in case (iii), $a_0 = b_0 > 0$ and $a_1 > 0$, forcing $b_1 = a_0 + a_1 > 0$ and then, regardless of the sign of $a_p$, no sign changes are introduced or removed in the passage from $a_0 + a_1, \ldots, a_p$ to $a_0, a_1, \ldots, a_p$ or from $b_1, \ldots, b_q$ to $b_0, b_1, \ldots, b_q$;
while in case (vi), $a_0 = b_0 > 0$,  $a_1 < 0$, $b_1 = a_0 + a_1 > 0$ and $a_p > 0$, so that in passing from $a_0 + a_1, \ldots, a_p$ to $a_0, a_1, \ldots, a_p$ two sign changes are introduced.

In each case, $\alpha - \beta$ is $0$ or $2$ and the lemma and the theorem are proved.

\subsection*{Remarks}

Several simple proofs of Descartes' rule of signs have been given over the years \cite{Gauss28,Albert43,Wang04,Komornik06}.
Our proof minimises the appeal to results such as the intermediate value theorem or Rolle's theorem.
In fact, it just uses the intermediate value theorem to show that $\PZ{f}=0$ if $\SC{f}=0$ and the rest is pure algebra.

The proof of the lemma works for an arbitrary ordered field.
This can be deduced from Descartes' rule of signs as formulated for an arbitrary real closed field by appeal to the existence of a real closure of an ordered field.
We have shown here that the rather heavy-weight machinery of real closed fields is not needed for this simple result.
\bibliographystyle{plain}
\bibliography{69}

\end{document}